\documentclass[11pt]{amsart} 

\usepackage{amsmath, amssymb, amsfonts, amsbsy, amsthm, array, latexsym, stmaryrd, fullpage, enumerate, ytableau, shuffle, multirow, mathdots, enumerate, tikz,ulem,youngtab,tikz,graphicx}
\usepackage{amsthm}
\usepackage{comment}
\usepackage{tikz}
\usetikzlibrary{matrix,arrows,decorations.pathmorphing}
\newtheorem{theorem}{Theorem}[section]

\theoremstyle{definition}

\newtheorem*{definition*}{Definition}
\newtheorem*{theorem*}{Theorem}
\newtheorem*{proposition*}{Proposition}
\newtheorem*{corollary*}{Corollary}

\begin{document}

\title{What is Schur positivity and how common is it?}
\author{Rebecca Patrias}
\maketitle

Peek into a lecture on algebraic combinatorics, representation theory, algebraic geometry, or theoretical physics these days and there is a good chance that you will hear something related to Schur positivity. This typically has to do with some interesting decomposition of an algebraic structure into irreducible components, which are naturally encoded by Schur polynomials. 

The Schur polynomials $s_\mu(x_1,x_2,\ldots, x_n)$ are symmetric functions that are indexed by weakly decreasing sequences of integers $\mu_1\geq \mu_2\geq \cdots\geq \mu_t>0$. Being symmetric polynomials means that they are invariant under all permutations of variables $x_i\mapsto x_{\sigma(i)}$. For instance, with $\mu=(2,1)$,
\begin{eqnarray*}
  s_{(2,1)}(x_1,x_2,x_3)&=&x_1^2x_2+x_1^2x_3+x_1x_2^2+x_1x_3^2+x_2^2x_3+x_2x_3^2+2\,x_1x_2x_3.
 \end{eqnarray*}
Non-zero Schur polynomials constitute a linear basis of the ring of symmetric polynomials. 
 
Named after Issai Schur, Schur polynomials originally occurred  in Jacobi's work as quotients of determinants:
 $$s_\mu( x_1,x_2,\ldots, x_n):=\frac{\det(x_i^{\mu_j+n-j})_{1\leq i,j\leq n}}{\det(x_i^{n-j})_{1\leq i,j\leq n}}.$$
This is indeed a polynomial, since both the numerator and denominator are divisible by all differences $(x_i-x_j)$, for $1\leq i<j\leq n$. The denominator is, in fact, equal to the product of all such differences and thus is the well-known Vandermonde determinant. Moreover, the symmetry of $s_\mu$ is obvious since the numerator and denominator are both antisymmetric.

It is not immediately clear from the definition above why Schur polynomials should play an important role in so many areas of mathematics. To help shed light on this question, we will see that Schur polynomials occur as characters of irreducible representations of the general linear group $GL(V)$. With this information in hand, it becomes apparent that their importance is derived from the fundamental role of $GL(V)$ in a wide spread of mathematical subjects. 
Let us now discuss the topics above more explicitly, with the goal of better understanding the central role of the Schur polynomials.

\section{Symmetric polynomials}

Consider the ring $\mathbb{Q}[x_1,\ldots,x_n]$ of polynomials in $n$ variables with rational coefficients. The symmetric group $S_n$ acts on this ring by permuting the variables: for $f(x_1,\ldots,x_n)\in\mathbb{Q}[x_1,\ldots,x_n]$ and $\sigma\in S_n$, \[ \sigma \cdot f(x_1,\ldots,x_n) = f(x_{\sigma(1)},x_{\sigma(2)},\ldots,x_{\sigma(n)}).\] We say that a polynomial $f$ is \textit{symmetric} if $\sigma\cdot f = f $ for all $\sigma\in S_n$. For instance, setting $n=3$, we have that  $x_1x_2+x_1x_3+x_2x_3$ is symmetric, while $x_1+x_2^2+x_3^2$ is not. 

There are many well-known $\mathbb{Z}$-bases for the ring of symmetric polynomials; let's build some intuition for the most natural of them. Say you are given a symmetric polynomial $f(x_1,x_2,x_3,x_4)$ and you know that the monomial $x_1^2x_2x_3$ is in $f$ with  coefficient $c$. Since $f$ is symmetric, we know that $x_1^2x_2x_4$ must also appear with coefficient $c$, and the same is true for $x_1^2x_3x_4$, $x_1x_2^2x_3$, $x_1x_2x_4^2$, etc. In other words, we know that any monomial of the form $x_1^{\alpha_1}x_2^{\alpha_2}x_3^{\alpha_3}x_4^{\alpha_4}$, where $\alpha=(\alpha_1,\alpha_2,\alpha_3,\alpha_4)$ gives $(2,1,1)$ when the nonzero components are sorted into decreasing order appears with coefficient $c$ in $f$. Let us call the sum of such monomials the \textit{monomial symmetric polynomial} $m_{(2,1,1)}(x_1,x_2,x_3,x_4)$. 

In general, let $\lambda=(\lambda_1,\ldots,\lambda_k)$ be a finite sequence of weakly decreasing positive integers. We call such a sequence a \textit{partition}. For $\lambda$ with at most $n$ parts, we can similarly define the monomial symmetric polynomial $m_\lambda(x_1,\ldots,x_n)$ by

\[m_\lambda(x_1,\ldots,x_n)=\sum_{\alpha}x_1^{\alpha_1}\cdots x_n^{\alpha_n},\]
where we sum over all distinct sequences of non-negative integers  $\alpha=(\alpha_1,\ldots,\alpha_n)$ that give $\lambda$ when the nonzero parts are arranged in decreasing order.  For example, \[m_{(4,2)}(x_1,x_2,x_3)=x_1^4x_2^2+x_1^4x_3^2+x_2^4x_3^2+x_1^2x_2^4+x_1^2x_3^4+x_2^2x_3^4\] 
since $(4,2,0)$, $(4,0,2)$, $(0,4,2)$, $(2,4,0)$, $(2,0,4)$, and $(0,2,4)$ are all the sequences with three  components whose nonzero parts arranged in decreasing order give $(4,2)$. 

After a little thought, we see that any monomial in our symmetric polynomial $f(x_1,x_2,x_3,x_4)$ belongs to some $m_\lambda(x_1,x_2,x_3,x_4)$ and if some monomial of $m_\lambda(x_1,x_2,x_3,x_4)$ appears in $f$, then every monomial of $m_\lambda(x_1,x_2,x_3,x_4)$ appears in $f$ with that same coefficient. It follows that we may decompose $f$ into a finite sum of monomial symmetric functions. In this way, we see that the set $\{m_\lambda(x_1,\ldots,x_n)\}$ forms a basis for the ring of symmetric polynomials in $n$ variables, where the collection includes all partitions with at most $n$ parts.

\section{Schur polynomial definition}
As mentioned before, the Schur polynomials form another basis for the ring of symmetric functions. Before going any further, let us take a brief detour into the life of their namesake, Issai Schur. Schur was born in Mogilev, Russia in 1875, and enrolled as a student at the University of Berlin in 1894. At the time, the University of Berlin was home to the well known Hermann Schwarz, Max Planck, and Ferdinand Frobenius, all of whom are mentioned at the end of Schur's dissertation. He became a professor at the Berlin academy, where his lectures were so popular that in the winter semester of 1930, the number of interested students exceeded the capacity of his 500-student lecture hall. However, not all of Germany felt so positively about Schur. As a Jewish mathematician, he was temporarily dismissed from his position in Berlin and was eventually forced to resign in 1935. His family relocated to Israel, where he continued his work. In 1940, Schur suffered a heart attack while giving a lecture, and he died weeks later. \cite{Schur}

There are many equivalent definitions of the Schur polynomials $s_\lambda$. Here we'll discuss a well known combinatorial definition as it will allow us to easily see the relation between Schur polynomials and monomial symmetric polynomials.

Associate to each partition  $\lambda$ a \textit{Young diagram}: a left-justified array of boxes where the $i$th row from the top has $\lambda_i$ boxes. Define a \textit{semistandard Young tableau} of shape $\lambda$ to be the result of filling each box of the Young diagram of shape $\lambda$ with a positive integer such that entries weakly increase from left to right along rows and strictly increase down columns. To each semistandard Young tableau, $T$, associate a monomial $x^T$, where the exponent of $x_i$ is the number of boxes labeled with $i$ in $T$. See Figure~\ref{fig:SSYT} for an example.

\begin{figure}[h]
\ytableausetup{boxsize=.2in}
$T=$
\begin{ytableau} 1 & 1 & 2 & 2 & 2 & 5 \\
2 & 4 & 4 & 4 \\ 4 & 5 \end{ytableau}$\hspace{.1in}\longleftrightarrow \hspace{.1in}x^T=x_1^2x_2^4x_4^4x_5^2$
\caption{A semistandard Young tableau of shape $(6,4,2)$ and the corresponding monomial.}\label{fig:SSYT}
\end{figure}

We then have have that 
\[s_\lambda(x_1,\ldots,x_n)=\displaystyle\sum_{\text{shape}(T)=\lambda} x^T,\] where we sum over all semistandard Young tableaux of shape $\lambda$ with entries in $\{1,2,\ldots,n\}$. For example, we see that $s_{(3,1)}(x_1,x_2)=x_1^3x_2+x_1^2x_2^2+x_1x_2^3$ by listing the semistandard Young tableaux with entries in $\{1,2\}$ below.

\begin{center}
\begin{ytableau}
1 & 1 & 1 \\ 2
\end{ytableau}\hspace{.3in}
\begin{ytableau}
1 & 1 & 2 \\ 2
\end{ytableau}\hspace{.3in}
\begin{ytableau}
1 & 2 & 2 \\ 2
\end{ytableau}\hspace{.3in}
\end{center}

The downside of this beautiful combinatorial definition is that, in contrast with the original definition as a quotient of determinants, it is not immediately clear why a Schur polynomial is symmetric.  One can use the Bender--Knuth involution, which describes how to exchange $x_i$ and $x_{i+1}$, to prove this symmetry using the tableau definition.

Given that the Schur polynomial $s_\lambda(x_1,\ldots,x_n)$ is symmetric, we should be able to write it in the basis of monomial symmetric polynomials. Suppose, for example, we would like to find the coefficient of $m_{(2,2,1)}(x_1,x_2,x_3)$ in $s_{(3,2)}(x_1,x_2,x_3)$. Using symmetry, it suffices to find  the coefficient of $x_1^2x_2^2x_3$ in $s_{(3,2)}(x_1,x_2,x_3)$. Using the tableau definition of the Schur polynomials, this amounts to counting the number of semistandard  Young tableaux of shape $(3,2)$ filled with two 1's, two  2's, and one 3. We call this number the \textit{Kostka number} $K_{(3,2),(2,1,1)}$, where in general, the Kostka  number $K_{\lambda,\mu}$ counts the number of semistandard tableaux of shape $\lambda$ filled with $\mu_i$ $i$'s. For any partition $\lambda=(\lambda_1,\ldots,\lambda_k)$, let $|\lambda|=\sum_{i=1}^k\lambda_i$. When $|\lambda|=r$, we say that $\lambda$  is a partition of $r$ and write $\lambda\vdash r$. Notice that $K_{\lambda,\mu}=0$ unless $|\lambda|=|\mu|$ because any Young diagram with $|\lambda|$ boxes needs exactly $|\lambda|$ entries to form a semistandard Young tableau. This idea leads us to the following.

\begin{theorem}
For any partition $\lambda$, \[s_\lambda(x_1,\ldots,x_n)=\displaystyle\sum_{\mu\vdash |\lambda|}K_{\lambda,\mu}m_\mu.\].
\end{theorem}

As we will see, these numbers play a key roll in understanding how likely it is for a symmetric polynomial to be Schur positive, i.e., to have an expansion in the Schur polynomial basis with positive coefficients.

\section{Polynomial representations of $GL(V)$}
One important place that Schur polynomials appear in nature is as the characters of irreducible representations of the complex general linear group. Let's see how this works.

Let $M,N\in\mathbb{N}$, $V$ be an $N$-dimensional complex vector space and $W$ be an $M$-dimensional complex vector space. A (finite-dimensional) \textit{representation} of the general linear group $GL(V)$ is a homomorphism $\varphi:GL(V)\to GL(W)$. In other words, $\varphi(AB)=\varphi(A)\varphi(B)$ and $\varphi(A)(w_1+w_2)=\varphi(A)w_1+\varphi(A)w_2$ for all $A,B\in GL(V)$ and $w_1,w_2\in W$. We may also choose bases for $V$ and $W$ and consider instead a homomorphism $\varphi: GL_N(\mathbb{C})\to GL_M(\mathbb{C})$. We then consider two representations to be equivalent if they differ by a change of basis. We will work in this latter setting so that we may use the language of matrices. 

We say that such a representation is \textit{polynomial} if the entries of the matrix $\varphi(A)$ are polynomials in the entries of $A$ for all $A\in GL_N(\mathbb{C})$. For example, the following $GL_2(\mathbb{C})$ representation is clearly polynomial. We leave it to the reader to check that it is indeed a homomorphism.

\[S^2: \begin{pmatrix}a & b \\ c & d\end{pmatrix}\longmapsto \begin{pmatrix} a^2 & 2ab & b^2 \\ ac & ad+bc & bd \\ c^2 & 2cd & d^2\end{pmatrix}\] 


Instead of focusing on a representation $\varphi$, we will explore its \textit{character}, char$_\varphi:GL_N(\mathbb{C})\to\mathbb{C}$, defined by char$_\varphi(A)=\text{trace}(\varphi(A))$. It turns out that char$_\varphi(A)$ depends only on the eigenvalues of the matrix $A$ and that the character of a polynomial representation will be a polynomial in the eigenvalues of the matrix $A$. We may thus consider char$_\varphi$ to be a polynomial char$_\varphi(x_1,\ldots,x_n)$. We can easily compute, for example, that if $A\in GL_2(\mathbb{C})$ has eigenvalues $\theta_1,\theta_2$, then $S^2(A)$ has eigenvalues $\theta_1^2,\theta_1+\theta_2,$ and $\theta_2^2$, and so char$_{S^2}(x_1,x_2)=x_1^2+x_1x_2+x_2^2$. We can now see that this character is a Schur polynomial: char$_{S^2}(x_1,x_2)=s_{(2)}(x_1,x_2)$.

Given two representations $\varphi_1$ and $\varphi_2$, we form their direct sum $\varphi_1\oplus \varphi_2$ by 
\[\varphi_1\oplus \varphi_2 (A):=\begin{pmatrix}\varphi_1(A) & 0 \\ 0 & \varphi_2(A) \end{pmatrix}\] for any $A$ in $GL_N(\mathbb{C})$. It follows immediately that char$_{\varphi_1\oplus\varphi_2}(x)=$char$_{\varphi_1}(x)+$char$_{\varphi_2}(x)$. We say that a representation $\varphi$ is \textit{irreducible} if $\varphi=\varphi_1\oplus \varphi_2$ implies $\varphi_1=0$ or $\varphi_2=0$. With this definition in hand, we can now discuss a connection between representation theory of $GL_N(\mathbb{C})$ and symmetric functions: the Schur polynomials in $N$ variables are precisely the characters of the irreducible polynomial representations of $GL_N(\mathbb{C})$.

\begin{theorem}
Irreducible polynomial representations $\varphi^\lambda$ of $GL_N(\mathbb{C})$ can be indexed by partitions $\lambda$ with at most $N$ components so that 
\[\text{char}_{\varphi^\lambda}(x_1,\ldots,x_N)=s_\lambda(x_1,\ldots,x_N).\]
\end{theorem}

From this theorem, it follows that Schur positivity is preserved by operations that correspond to operations on representations such as addition and multiplication, as well as more involved operations like plethsym and the Kronecker product, which we have not defined here. With lots of ways to produce new Schur positive expressions, one might guess that Schur positivity must be pretty common.

\section{Probability of Schur positivity}
By this point, we have seen that Schur positivity is a natural phenomenon that is worthy of study. We now seek to answer the question of the prevalence of Schur positive symmetric polynomials. Are most symmetric polynomials Schur positive? Or are proofs of Schur positivity of interest because the property is very rare? It turns out that the latter is the case. Let's see what we mean by this.

Suppose we are handed a homogeneous symmetric polynomial of degree $k$ with nonnegative coefficients $f(x_1,\ldots,x_n)$. (Any negative coefficients would of course immediately rule out the possibility of Schur positivity.) We would like to answer the following question: What is the probability that $f(x_1,\ldots,x_n)$ is Schur positive?

Let us first make this question more precise. Any such polynomial may be expressed as a nonnegative linear combination of monomial symmetric polynomials indexed by partitions of $k$
\[f(x_1,\ldots,x_n)=\displaystyle\sum_{\lambda\vdash k} a_\lambda m_\lambda,\] and the set of such polynomials forms a cone $$M=\left\{\displaystyle\sum_{\lambda\vdash k}  a_\lambda m_\lambda(x_1,\ldots,x_n) \mid a_\lambda\in\mathbb{R}_{\geq0}\right\}.$$ Inside the cone $M$ sits the cone of Schur positive homogeneous symmetric polynomials of degree $k$:
$$S=\left\{\displaystyle\sum_{\lambda\vdash k}  a_\lambda m_\lambda(x_1,\ldots,x_n) \mid a_\lambda\geq 0\text{ and $f$ is Schur positive}\right\}.$$

The answer to our question can thus be found by calculating the ratio of the volume of $S$ to the volume of $M$. We can do this by realizing that the ratio of the volumes will be the same as the ratio of the respective volumes of any slice of $M$.

To this end, we consider the slice consisting of the functions whose coefficients add to one when expressed in the monomial basis
$$M_{\text{slice}}=\left\{\displaystyle\sum_{\lambda\vdash k}  a_\lambda m_\lambda(x_1,\ldots,x_n) \mid a_\lambda\in\mathbb{R}_\geq 0 \text{ and } \sum a_\lambda =1 \right\}$$ and the corresponding slice $S_{\text{slice}}$. For the remainder of this explanation, we restrict ourselves to homogeneous symmetric polynomials of degree $k=3$ for simplicity of explanation. The argument follows the same way for arbitrary $k$.

We first wish to find the volume---or area in this case---of $M_{\text{slice}}$. Let's fix the point $m_{(1,1,1)}$ to be the origin of $M_{\text{slice}}$ and consider the vectors 
\begin{eqnarray*} e_{1}&:=&m_{(2,1)}-m_{(1,1,1)}  \\ e_{2}&:=&m_{(3)}-m_{(1,1,1)} .\end{eqnarray*} 

It is clear that $\{e_1,e_2\}$ forms a basis for $M_{\text{slice}}$ since the monomial symmetric functions are linearly independent. Then the volume of $M_{\text{slice}}$ is easily seen to be 
\[\frac{1}{2}|\det(e_1,e_2)|=\frac{1}{2}.\]

Now consider some $s_\lambda$ for $\lambda\vdash 3$, and define $k_\lambda = \sum_\mu K_{\lambda\mu}$. We then have $k_{(3)}=3$, $k_{(2,1)}=3$, and $k_{(1,1,1)}=1$, where, for example, $k_{21}=3$ 
since $s_{21}=m_{21}+2m_{111}$. 
Then the points \[\left\{\frac{1}{k_\lambda}s_\lambda\mid \lambda\vdash n\right\}=\left\{\frac{1}{3}s_{(3)},\frac{1}{3}s_{(2,1)},s_{(1,1,1)}\right\}\] are the vertices that determine $S_{\text{slice}}$. Noticing that $s_{(1,1,1)}=m_{(1,1,1)}$, we have that $s_{(1,1,1)}$ is the origin of $S_{\text{slice}}$. We define basis vectors 
\begin{eqnarray*}v_1&=&\frac{1}{3}s_{21}-s_{111}=\frac{1}{3}m_{21}+\frac{2}{3}m_{111}-m_{111} =\frac{1}{3}e_1 \\ v_2&=&\frac{1}{3}s_{3}-s_{111}=\frac{1}{3}m_{3}+\frac{1}{3}m_{21}+\frac{1}{3}m_{111}-m_{111} = \frac{1}{3}e_1+\frac{1}{3}e_2.\end{eqnarray*}

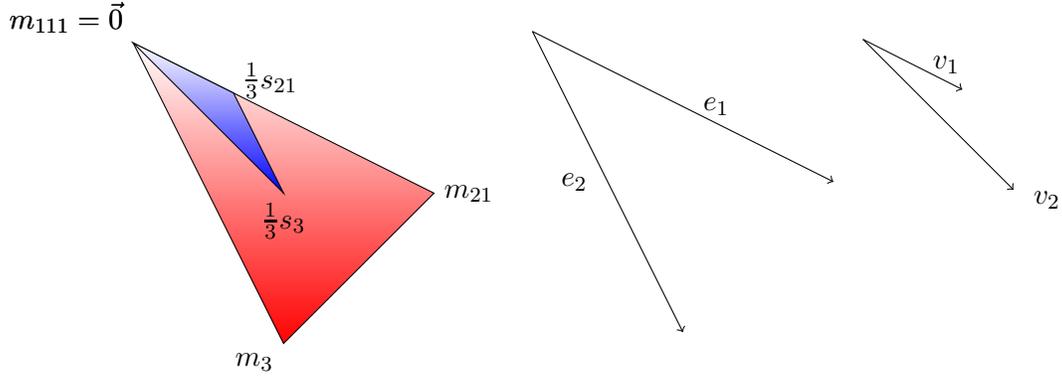
\begin{figure}
\begin{center}
\begin{tikzpicture}[scale=2]
\draw[top color=white,bottom color=red] (0,0) -- (2,-1) -- (1,-2) -- cycle; 
\foreach \coordinate/\label/\pos in {{(0,0)/$m_{111}=\vec{0}$/above left},{(2,-1)/$m_{21}$/right},{(1,-2)/$m_3$/below left}}
\node[\pos] at \coordinate {\label};
\draw[top color=white,bottom color=blue] (0,0) -- (2/3,-1/3) -- (1,-1) -- cycle; 
\foreach \coordinate/\label/\pos in {{(0,0)/$m_{111}=\vec{0}$/above left},{(.66,-.25)/$\frac{1}{3}s_{21}$/right},{(1,-1)/$\frac{1}{3}s_3$/below}}
\node[\pos] at \coordinate {\label};
\end{tikzpicture}\hspace{.1in}
\raisebox{.25in}{\begin{tikzpicture}[scale=2] 
\draw[->] (0,0) -- (2,-1)  node [midway, label=right:$e_1$] {};
\draw[->] (0,0) -- (1,-2) node [midway, label=left:$e_2$] {};
\end{tikzpicture}}\hspace{.1in}
\raisebox{.85in}{\begin{tikzpicture}[scale=2] 
\draw[->] (0,0) -- (1,-1) node [below, label=right:$v_2$] {};
\draw[->] (0,0) -- (.66,-1/3) node [midway, label=right:$v_1$] {};
\end{tikzpicture}}
\end{center}
\caption{An illustration of $S_{\text{slice}}$ inside of $M_{\text{slice}}$ in degree 3. }
\label{fig:vectors}
\end{figure}
Therefore the volume (area in this case) of $S_{\text{slice}}$ is
\[\frac{1}{2}|\det(v_1,v_2)|=\left(\frac{1}{2}\right)\left(\frac{1}{3}\right)\left(\frac{1}{3}\right)=\frac{1}{2}\prod_{\lambda\vdash 3} \frac{1}{k_\lambda}.\]
Hence the ratio of the volume of the slice is $S_{\text{slice}}/M_{\text{slice}}=(\frac{1}{3})(\frac{1}{3})=\prod_{\lambda\vdash 3} \frac{1}{k_\lambda}$. In general, we have the following.

\begin{theorem}[F. Bergeron, R. Patrias, and V. Reiner] Let $f$ be a homogeneous symmetric function of degree $k$ with non-negative coefficients. Then the probability that $f$ is Schur positive is \[\prod_{\lambda\vdash k} \frac{1}{k_\lambda}.\]
\end{theorem}

As $k$ grows, we can now see that the probability that $f$ is Schur positive rapidly becomes quite small. Indeed, for $k=1$ through $k=7$, the probabilities are 1, $\frac{1}{2}$, $\frac{1}{ 9}$, $\frac{1}{560}$, $\frac{1}{480480}$, $\frac{1}{1027458432000}$, and $\frac{1}{2465474364698304960000}$.

This is a typical phenomenon: a rare notion becoming prevalent in the presence of interesting structure. Next time you see a Schur positive symmetric polynomial, you'll know it is a rare beast with representation theory lurking behind it.

\end{document}